
\magnification=1200 
\let\mathcal=\cal
\def\sqr#1#2{{\vcenter{\vbox{\hrule height.#2pt
\hbox{\vrule width.#2pt height #1pt \kern#1pt
\vrule width.#2pt}
\hrule height.#2pt}}}}
\def\square{\mathchoice\sqr56\sqr56\sqr{2.1}3\sqr{1.5}3}
\def\qed{\hfill$\square$}
\def\IR{\hbox{\rm I\kern-.2em\hbox{\rm R}}}

\vskip1.4cm \noindent

\centerline  {\bf Geometric lower bounds for the spectrum 
of elliptic PDEs}
\par
\centerline {\bf with Dirichlet 
conditions in part}
\bigskip
\par
\centerline {{ \bf Evans M. Harrell II}\footnote{}{harrell@math.gatech.edu, 
School of Mathematics,Georgia Tech, Atlanta, GA 30332-0160, USA.  
\copyright 2004 by the author. Reproduction of this article, in its 
entirety, by any means is permitted for non--commercial purposes.  
This work was supported by NSF grant DMS-0204059.}}   \smallskip

\vskip .7 true in

 \centerline{\bf Abstract}

An extension of the lower-bound lemma of Boggio is 
given for the weak forms of certain elliptic operators,
which have partially Dirichlet and partially Neumann
boundary conditions, and are in general nonlinear.  Its consequences and those of an adapted Hardy inequality for the location of the bottom of the spectrum are 
explored in corollaries
wherein a variety of assumptions 
are placed on the shape of the Dirichlet and Neumann boundaries.
\vskip .7 true in

\vfill\eject 
\noindent{\bf I. Introduction} 
\medskip E.B. Davies has promoted the use of operator and quadratic--form 
inequalities to obtain 
spectral information about Laplacians and Schr\"odinger operators.  Among the goals of his analysis 
has been to understand the effect of the shape of the boundary on eigenvalues and the
eigenfunctions.
The present note may be viewed as a
scholium to parts of the first chapter of [Dav89], where 
linear elliptic operators are bounded below, in the weak sense, by functions 
related to the shape of the boundary of a domain.
The aim is the elaboration of lower bounds of
similar kinds, given two complicating features:  boundary conditions other 
than of Dirichlet type, and nonlinearities in the principal part of the operators.  

As has been pointed out in [FHT99], a precursor to many inequalities used for lower bounds in the spectral theory of elliptic differential equations is to be found a 1907 article by T. Boggio [Bog07].  This direct
consequence of the Gau\ss--Green theorem
was interpreted in [FHT99] as a variational principle for the Laplacian and its nonlinear 
generalisation the$p$--Laplacian, on bounded domains with vanishing Dirichlet conditions on the boundary.  
Lemma I.1 below, the main tool used in this article, likewise applies to
a family of nonlinear elliptic operators, and non--Dirichlet boundary conditions are allowed.
No great originality is claimed for Lemma I.1:  The idea is still
close to that of
[Bog07], and it is likely that an equivalent result could be extracted from what 
has been reported in
[Maz81] or [OpKu90].  It is, however, stated in a form allowing geometric spectral results to be
extracted easily, with attention to 
non--Dirichlet boundary portions.  The short proof is given so that the exposition is self--contained.

The great majority of lower bounds to the spectrum of elliptic partial differential equations suppose Dirichlet conditions throughout any boundary that is present.  One reason for this is that if
Neumann conditions occur, the spectrum of the Laplacian 
is highly unstable with respect to perturbations of the boundary.  An arbitrarily small perturbation of the Neumann part of the boundary may cause the lowest eigenvalue to decrease arbitrarily 
close to $0$.  Intuitively, this is possible because the lowest eigenvalue of any domain with pure Neumann conditions is 0, so if a small such domain is weakly coupled to $\Omega$ via a thin neck or partitions, then
the lowest eigenvalue will remain tiny.  
The instability of Neumann spectra 
has been extensively investigated, and instructive ``rooms and passages''
models date from Courant and Hilbert:  See [CoHi37], [Maz73], {[Fra79], [EvHa89],[HSS91],[Arr95], [BuDa02] among many treatments of the Neumann spectrum.

There is evidence that if assumptions on the
Neumann part of the boundary prevent its near--isolation, then the effect on the lowest eigenvalue
will be more moderate.
For instance, in [BuDa02], a 
continuous perturbation theory is established for purely Neumann Laplacians with respect to 
uniformly H\"older perturbations.  Section II of this article includes lower bounds controlling the downward shift in the fundamental eigenvalue upon 
non--perturbative
enlargement of the domain.

It will be assumed throughout that $\Omega$ is a bounded 
domain in $\IR^d$, the boundary of which is sufficiently
regular that Green's formula is valid for suitable test
functions and that
the outward normal vector,
denoted ${\bf \nu}$, exists a.e. at any part of the boundary where conditions other than 
vanishing Dirichlet
are imposed.  Regularity depends on the validity of certain Sobolev 
embeddings.  
Here each domain $\Omega$ will be assumed to have a 
boundary consisting of two pieces each of positive $d-1$--dimensional Lebesgue measure, denoted
$\partial \Omega_D$ and $\partial \Omega_N := \partial \Omega \setminus \partial \Omega_N$. 
For the purposes of this article it suffices
to suppose that the boundary is of ``class $C$'' (see [EdEv87], p. 244, 248), and that the boundary 
is $C^1$ a.e. with respect to the surface Lebesgue measure.
Such a domain will be termed {\it regular}.   (Regularity is not
actually required on $\partial \Omega_D$, while on 
$\partial \Omega_N$ usefully wider conditions 
of regularity may possibly be gleaned from [Ken94].  The widest
conditions for the boundary shall not be pursued here.)  
An adapted set of test functions ${\mathcal D}\left({\Omega, \Omega_D}\right)$
can be defined as consisting of the restrictions to $\Omega$ of 
$C_c^{\infty}\left({\IR^d  \setminus \partial\Omega_D}\right)$.

Although the notation $\partial \Omega_N$ is introduced with Neumann boundary conditions in 
mind, no such assumption is made.  
Only the Dirichlet portion of the boundary is actually 
subject to (weakly defined) boundary conditions.   
This article is concerned solely with functionals
on ${\mathcal D}$, without investigating the domains of
definition of the elliptic operators that the functionals may define.

\proclaim Lemma I.1.
Suppose that $\Omega$ is regular, 
$1 < p < \infty$, 
${\bf Q}$ is a real--valued 
absolutely continuous vector field, and $A$ is a 
real--valued absolutely continuous tensor field
of type ${\mathcal T}_1^1$.  {\rm (In elementary terms,
$A$ is a $d \times d$ real matrix--valued function acting on 
${\bf Q}$.)}  Suppose further that
${\bf \nu} \cdot A^T {\bf Q} \leq 0$ a.e. on $\Omega_N$.
Then for all $\zeta \in {\mathcal D}\left({\Omega, \Omega_D}\right)$,
$$\int_{\Omega }^{}{{\left|{A \nabla \zeta }\right|}^{p}dx}\ \ge  \ \int_{\Omega }^{}\left\{{{\rm div}\left({A^T \bf Q}\right)\  -\ (p-1)  {\left|{\bf Q}\right|}^{p^{\prime}}}\right\}{\left|{\zeta }\right|}^{p}\ dx,\eqno{(1.1)}$$
where $p'$ is the dual index, $p' = {{p}\over {p-1}}$.  
\par\noindent
{\it
Alternatively, Let $a$ and $q$ be absolutely continuous functions on $\Omega$ 
and let ${\bf u}$ be a fixed unit vector
such that 
$a q {\bf u} \cdot {\bf \nu} \le 0$ a.e. on $\partial \Omega_N$.
Then for all $\zeta \in {\mathcal D}\left({\Omega, \Omega_D}\right)$,}
$$\int_{\Omega}^{}{\left|a {{\bf u} \cdot \nabla \zeta}\right|^p dx} \geq \int{\left\{{{{\bf u} \cdot \nabla(a q)} - (p-1) q^{p'}}\right\} |\zeta|^p dx}.\eqno{(1.2)}$$

\noindent {\bf Proof}.
We calculate following [FHaT99] with the divergence theorem:
$$\eqalign{
0 &\leq - \int_{\partial \Omega}{ {\bf \nu} \cdot 
\left(A^T {\bf Q}|\zeta|^p\right) dS}
= -\int_{\Omega} {{\rm div} \left(A^T {\bf Q}|\zeta|^p\right) dx}\cr
&=-\int_{\Omega}{({\rm div}  A^T {\bf Q}) |\zeta|^p dx}-
p\int_{\Omega}{|\zeta|^{p-2}\zeta{\bf Q} \cdot A \nabla \zeta dx.}}$$ 
With ${\bf w}=|\zeta|^{p-2}\zeta\bf Q$, the arithmetic--geometric mean inequality 
[HLP34] gives 
$$|\left(A \nabla \zeta\right) \cdot {\bf w}|\le {1\over p}\, 
|A \nabla \zeta|^p+{1\over p'}\, |{\bf w}|^{p'}={1\over p}\, 
|A \nabla \zeta|^p +{1\over p'}\,|\zeta|^p |{\bf Q}|^{p'}.$$ 
Collecting terms,
$$\int_\Omega{|A \nabla \zeta|^p dx} \ge \int_\Omega 
{\left({\rm div} \left({A^T {\bf Q}}\right)-\left({p-1}\right) \, 
|{\bf Q}|^{p'}\right)|\zeta|^p dx},$$
which is (1.1).  The proof of (1.2) is strictly analogous, beginning with
the Gau\ss--Green formula in the version
$$0 \le - \int_{\partial \Omega} {\left({a q |\zeta|^p}\right) \nu^i dS}= - \int_{\partial \Omega} {{{\partial\left(a q |\zeta|^p \right)} \over {\partial x_i}}dx}$$
for a Euclidean coordinate ${x_i}$ (e.g., [Ev98]).
The result is
$$\int_{\Omega}^{}{\left|a {{\partial \zeta} \over {\partial x_i}}\right|^p dx} \geq \int{\left\{{{{\partial (a q)} \over {\partial x_i}} - (p-1) q^{p'}}\right\} |\zeta|^p dx},$$
which when written without coordinates becomes (1.2).
\hfill\qed
\smallskip
\par\noindent
{\bf Remarks}
{\parindent=18pt
\item{1.}
Inequality (1.1) is well--known in the case $p=2$ with Dirichlet boundary conditions on the
entire boundary, in the form obtained with 
$${\bf Q}_{std} := -A {\bf \nabla} \log(\Phi),$$ 
for a suitable function 
$\Phi > 0$ on $\Omega$.  
For $p \ne 2$ the analogous choice is
$${\bf Q}_{std} = - {{\left|{A \nabla \Phi}\right|^{p-2} A \nabla \Phi} \over {\Phi^{p-1}}}.$$
A simple calculation shows that for this choice
Lemma I.1 becomes 
$$\int_{\Omega }{{\left|{A \nabla \zeta }\right|}^{p} dx}\ \ge \  \int_{\Omega }{{-{\rm div} \left({\left|{A \nabla \Phi}\right|^{p-2} 
A^T A \nabla \Phi}\right) \over {\Phi^{p-1}} \left|{\zeta}\right|^p} dx},
\eqno(1.3)$$
which for $p=2$ is essentially Theorem 1.5.12 of  [Dav89].  
Let us refer to this as the {\it standard form} of Inequality (1.1).  It is all the more familiar
for the linear Dirichlet Laplacian, i.e., with $A$ taken as the identity and $p=2$, in which case
it becomes the well-known bound
$$- \Delta \ge {{- \Delta \Phi} \over {\Phi}},$$
which is often attributed to Barta or Duffin, who, 
however, published much later than Boggio [Bog07].
(In the not widely accessible
[Bog07], Boggio used the Gau\ss-Green Theorem to 
prove the case of Lemma I.1 with $d=2=p$, pure Dirichlet, $A = $ the identity.)
Lemma I.1 is truly more general than the 
standard form, which corresponds to a vector field ${\bf Q}$ restricted to be A times something
irrotational.
\item{2.}
Similarly, if 
$$q_{std} := {{-\left|a {{\partial \Phi}\over{\partial x_i}}\right|^{p-2} {{\partial \Phi}\over{\partial x_i}}} \over {\Phi^{p-1}}}$$
is inserted into (1.2), the result is
$$\int_{\Omega}^{}{\left|a {{\partial \zeta} \over {\partial x_i}}\right|^p dx} \geq \int_{\Omega}{{{-{{\partial}\over{\partial x_i}}\left({\left|a {{\partial \Phi}\over{\partial x_i}}\right|^{p-2} {{\partial \Phi}\over{\partial x_i}}}\right)} \over {\Phi^{p-1}}} |\zeta|^p dx},\eqno(1.4)$$
provided that $a q \nu^i \le 0$ a.e. on $\partial \Omega_N$.  The expression in the numerator of the
lower bound (1.4) is a one-dimensional $p$--Laplacian.  
\item{3.}
Lemma I.1 may be thought of as a  lower--bound variational principle for the 
bottom of the spectrum
of the operators whose energy forms are
as in the left side of (1.1), i.e.,
	$$H \zeta := - {\rm div} \left({\left|{A \nabla \zeta}\right|^{p-2} A^T A \nabla \zeta}\right).\eqno(1.5)$$
If
$\Phi > 0$ is the fundamental eigenfunction of $H$, 
that is, $H \Phi = \lambda_1 \Phi^{p-1}$,
then the lower bound to the spectrum of $H$ thus obtained reduces to
$\lambda_1$.
In this sense the inequality is optimal.   
While it must be conceded that the inequality is 
thus arguably not better in theory than the standard form
for the purpose of finding lower bounds, it has two practical advantages:
{\parindent 36 pt
\item{a)}
It is relatively easy to construct a vector field with positive divergence
to cover a domain the only
material restriction being that it be incoming on $\partial \Omega_N$.
\item{b)}
The nonlinearity of the bound can be exploited.  If {\it any} vector field can be 
constructed on
$\Omega$ with the necessary boundary behavior and
strictly positive divergence, then by scaling ${\bf Q} \rightarrow t {\bf Q}$ a strictly positive 
lower bound will be obtained for sufficiently small $t > 0$.

}
\item{4.}
An extension from Euclidean domains to subdomains of orientable smooth 
manifolds is straightforward.

}
\medskip

In the case where $A$ is the identity tensor, $H$ corresponds to the $p$--Laplacian, denoted 
$- \Delta_p$.
It is known
that the minimum of the functional 
$$\int_{\Omega }^{}{\left|{ \nabla \zeta }\right|^{p} dx},$$
\noindent
for $\|\zeta\|_p = 1$ and Dirichlet boundary conditions, corresponds to an eigenvalue $\lambda_1$ of
$- \Delta_p$
Many other familiar facts for the linear case $p=2$ carry over, such as that the minimiser $u_1$
exists in the appropriate Sobolev space and is positive on
$\Omega$. 
The second eigenvalue of (1.5) is also inflection point of (1.5) and can thus
be characterised variationally.  However, the theory of the 
higher eigenvalues of $- \Delta_p$ remains murky when $p \ne 2$.   
Detailed 
analysis of the $p$-Laplacian will not be entered into here, as the focus will be entirely on
lower bounds for its energy form $H$.
(The spectral theory of linear elliptic differential operators is presented, for example,
in [EdEv87], [Dav95], and a useful review of the
spectral theory of the $p$--Laplacian is to be found in [DKN96, DKT99].)

In the following section particular vector fields 
${\bf Q}$ will be chosen responding to assumptions about the shape of the 
domain $\Omega$.
\bigbreak

\noindent{\bf II. Spectral bounds with Dirichlet conditions on subsets of the boundary} 
\medskip

The aim of this section is to find lower bounds to the spectrum in terms of the shapes of 
$\partial \Omega_D$ and $\partial \Omega_N$.
For simplicity it will be assumed henceforth that $A$ is the identity tensor,
i.e., the case of the $p$--Laplacian.  No material difficulty would
arise if $A$ were
retained.  

In essence this section consists of a selection of examples of the use of Lemma I.1,
organized as a list of corollaries.

\medskip
\par\noindent
{\bf Definition II.1}.  
Let $\lambda_1\left({\Omega, \partial \Omega_D}\right)$ denote the infimum for
$\zeta \in {\mathcal D}\left({\Omega, \partial \Omega_D}\right)$ (and not identically
zero) of
$$E_p(\zeta) := \left({{\int_{\Omega }^{}{\left|{\nabla \zeta }\right|}^{p}} \over {\int_{\Omega }^{}{\left|{ \zeta }\right|}^{p}}}\right).$$
The quantity $\lambda_1\left({\Omega, \partial \Omega_D}\right)$ will be referred to as the
{\it fundamental
eigenvalue
of the $p$--Laplacian}
with respect to $\Omega, \partial \Omega_D$.

As remarked in the introduction, eigenvalues can be rather unstable with 
respect to Neumann conditions on the boundary.  
In the purely
Dirichlet case the principle of domain--monotonicity is familiar:  If $\Omega \subset \Upsilon$,
then $\lambda_1\left({\Omega, {\rm pure DBC}}\right) \ge \lambda_1\left({\Upsilon, {\rm pure DBC}}\right)$.  
There is no such
general principle in the presence of Neumann conditions.  Nonetheless, Lemma I.1 implies the following:

\proclaim Corollary II.2 -- Restricted Monotonicity Principle.
Let $\Omega$ and $\Upsilon$ be regular domains such that
$\Omega \subset \Upsilon$.  Let $u_1(\Upsilon,\partial\Upsilon_D,p) > 0$ on $\Omega$
be a fundamental eigenfunction
for the $p$--Laplacian with respect to $\Upsilon, \Upsilon_D$.
If $\nu \cdot \nabla u_1(\Upsilon, \partial\Upsilon_D,p) \ge 0$ a.e. on $\partial \Omega_N$,
then $\lambda_1\left({\Omega, \partial\Omega_D, p}\right) \geq \lambda_1\left({\Upsilon, \partial\Upsilon_D, p}\right)$.

\par\noindent
{\bf Proof}.
The assumptions allow the standard form (1.3) of the lower bound
with the choice $\Phi = u_1(\Upsilon, \partial\Upsilon_D,p)$. 
In this case the right side of (1.3) 
reduces to $\lambda_1\left({\Upsilon, \partial\Upsilon_D, p}\right).$
\hfill \qed
\medskip
An immediate consequence is that if the Neumann boundary is a graph with respect 
to a hyperplane, there is an analogue of the 
elementary outradius bound for the Dirichlet problem.
There results a crude ``box bound'' in terms of $\mu_{I,p} :=$ the fundamental eigenvalue of the 
one--dimensional $p$--Laplacian $-\Delta_p^{(1)}$ on the unit interval, with Dirichlet conditions at one
end and Neumann at the other.  (Thus $\mu_{I,2} = {\pi^2 \over 4}$.)

\proclaim Corollary II.3 -- Box Bound.
Let $\Omega$ be a regular domain such that 
$\partial \Omega_N$ is the graph of a function  on a subset of $\left\{{\bf x}: x_i = 0\right\}$, with 
${\bf \nu} \cdot {\bf e}_i \le 0$ a.e. on $\partial \Omega_N$  For $L := sup\left({x_i : {\bf x} \in \Omega}\right) - inf\left({x_i : {\bf x} \in \Omega}\right)$, 
$$\lambda_1\left({\Omega, \partial \Omega_D,p}\right) \ge {{\mu_{I,p} }\over{L^p}}$$
\smallskip
\par\noindent
{\bf Proof}.
Since the fundamental eigenfunction $u_{1,p}(x)$ of the one--dimensional $p$--Laplacian with Dirichlet conditions at one end and Neumann at the other is positive and monotonic 
(due to the maximum principle), we may choose $a=1$, $\Phi = u_{1,p}$ in (1.4), yielding the result.
\hfill \qed

\medskip
Less crude bounds involving the shape of 
$\Omega$ can be obtained with more sophisticated comparisons.  Recall that
in [Dav98], Davies exploited Hardy's one--dimensional inequality 
[Har20] [HLP34] to produce
lower bounds in the quadratic--form sense of the type
$$-\Delta \ge {C \over {\left({m({\bf x})}\right)^2}},$$ 
where $m$ is an averaged distance to the boundary, supposed entirely Dirichlet.
These bounds extend almost immediately to the case of the $p$--Laplacian and boundaries
that are only partly Dirichlet, with certain restrictions:

\par\noindent
{\bf Definition II.4}.  
Given
a regular domain $\Omega$, the
{\it Dirichlet boundary sector} of ${\bf x}, \partial \Omega_D$,
denoted $S\left({\bf x}, \partial \Omega_D\right)$,
is the union of all line segments contained
in $\Omega$ joining ${\bf x}$ to $\partial \Omega_D$. 
Following [Dav89], an averaged distance to the Dirichlet boundary is given by 
$${1 \over {\left({m\left({\bf x}\right)}\right)^p}}
:= {1 \over \omega_d} \int_{K}{{dS({\bf u})} \over {\left({d_{\bf u}({\bf x}})\right)^p}}
$$
where the measure $dS({\bf u})$ is the Lebesgue measure 
on the unit sphere $S^{d-1}$, 
$\omega_d$ is the total measure of $S^{d-1}$,
$K := S\left({\bf x}, \partial \Omega_D\right) \cap S^{d-1}$,
and $d_{\bf u}$ is the length of a line segment joining ${\bf x}$ to $\partial \Omega_D$.
\par\smallskip
\noindent
Note:  The Dirichlet boundary sector of ${\bf x}$ may be vacuous, and $1 \over m$ may be zero.

\proclaim Lemma II.5.
Let $\Omega$ be a regular domain and $\zeta \in {\mathcal D}\left({\Omega, \Omega_D}\right)$.  
Then for any  $1 < p < \infty$,

$$\int_{\Omega}{\left|\nabla \zeta\right|^p dx} \ge p^{-p}(p-1)^{p-2} \int_{\Omega}{\left|{\zeta({\bf x})} \over {m({\bf x})}\right|^p  dx}.$$
\medskip
\par\noindent
With the one--dimensional $L^p$ Hardy inequality 
(i.e., [OpKu90] Theorem 1.14 with $p=q$, $v=1$, and $w(x) = (x-a)^{-p}$, itself derivable from 
lemma I.1), the
proof of the lemma is almost exactly as for [Dav98], Theorem 1.5.3, and will 
be omitted.

%

Suppose that $\partial \Omega_N$ is starlike with respect to an exterior point, at which the origin will be placed.  It will now be shown that there is a sort of
multidimensional Hardy inequality with respect to a distance function that is the intuitively the lesser of
a multiple of
$|{\bf x}|$ and the distance to the Dirichlet part of the boundary.  The precise version of this
statement is expressed by (2.2) in Corollary II.6 and (2.3) in Corollary II.7.  
Versions of these corollaries, in the linear case $p=2$, have been
previously reported, though not published, by the 
author since [Har93].  They imply, for example, that it is possible to have
a fractal $\partial\Omega_N$ without the collapse of the bottom of the spectrum to $0$.

\proclaim
Corollary II.6.
Let $1 < p < d$, and
suppose that the origin is exterior to $\Omega$, that 
$\partial \Omega_N$ is starlike with respect to the origin, and that the interior of $\Omega$ adjacent to 
$\partial \Omega_N$ is away from the origin.  Then for all
$\zeta \in {\mathcal D}\left({\Omega, \partial\Omega_D, p}\right)$,
$$\int_{\Omega }^{}{\left|{\nabla \zeta }\right|^{p} dx} 
\ge \left({{d-p}\over p}\right)^p \int_{\Omega }
{{\left|{\zeta \over \left|{\bf x}\right|}\right|}^{p} dx}.\eqno(2.1)
$$
\par\noindent
{\it Consequently,
$$\lambda_1\left({\Omega, \partial \Omega_D,p}\right) \ge \left({{d-p} \over {p \inf\left\{|{\bf x}| : {\bf x} \notin \Omega\right\}}}\right)^p.$$
Moreover, for any} $\gamma, 0 \le \gamma \le 1,$
$$\lambda_1\left({\Omega, \partial \Omega_D,p}\right) \ge \inf_{{\bf x} \in \Omega}{\left({\gamma \left({{d-p} \over {p |{\bf x}|}}\right)^p+ (1-\gamma) {(p-1)^{p-2} \over {{p^{p} m({\bf x})^p}}}}\right)}.\eqno(2.2)$$

\noindent
{\bf Proof}.
Eq. (2.1) results from the choice ${\bf Q} = \left({{d-p} \over p}\right)^{p-1} {{\bf x} \over {\left|{\bf x}\right|^p}}$ in Lemma I.1 and a calculation.  The next statement is the immediate lower bound 
by replacing a factor in the integrand by its infimum.  Statement (2.2) is obtained in the same way
from the
weighted average of the lower bounds of Lemma II.5 and (2.1).
\hfill\qed
\medskip
When $p > d$, the origin needs to be taken from the other side of
$\partial \Omega_N$:
\medskip

\proclaim
Corollary II.7.
Let $d < p < \infty$, and
suppose that the origin is exterior to $\Omega$, that 
$\partial \Omega_N$ is starlike with respect to the origin, and that the exterior of $\Omega$ adjacent to 
$\partial \Omega_N$ is away from the origin.  Then for all
$\zeta \in {\mathcal D}\left({\Omega, \partial\Omega_D, p}\right)$,
$$\int_{\Omega }^{}{\left|{\nabla \zeta }\right|^{p} dx} 
\ge \left({{p-d}\over p}\right)^p \int_{\Omega }
{\left|{\zeta \over \left|{\bf x}\right|}\right|^{p} dx}.\eqno(2.3)
$$
\par\noindent
{\it Consequently,
$$\lambda_1\left({\Omega, \partial \Omega_D,p}\right) \ge \left({{p-d} \over {p \inf\left\{|{\bf x}| : {\bf x} \notin \Omega\right\}}}\right)^p.$$
Moreover, for any} $\gamma, 0 \le \gamma \le 1,$
$$\lambda_1\left({\Omega, \partial \Omega_D,p}\right) \ge \inf_{{\bf x} \in \Omega}{\left({\gamma \left({{p-d} \over {p |{\bf x}|}}\right)^p+ (1-\gamma) {(p-1)^{p-2} \over {{p^{p} m({\bf x})^p}}}}\right)}.\eqno(2.4)$$
\medskip
\noindent
{\bf Proof}.
Eq. (2.3) results from the choice ${\bf Q} = -\left({{p-d} \over p}\right)^{p-1} {{\bf x} \over {\left|{\bf x}\right|^p}}$ in Lemma I.1 and a calculation.  The other statements follow as in Corollary II.6.
\hfill\qed
\bigskip
For the missing case $p=d$, a comparison can be made with 
annular domains, i.e., the region between two concentric circles or spheres.

\medskip
\par\noindent
{\bf Definition II.8}.
Let $\alpha(r,R,p,d)$ denote the lowest eigenvalue of the $p$--Laplacian on 
$A := \{{\bf x}: r < \left|{{\bf x}}\right) < R\}$, with Dirichlet
boundary conditions imposed where $\left|{{\bf x}}\right| = R$ and
Neumann
boundary conditions imposed where $\left|{{\bf x}}\right| = r$.  
Similarly, let $\beta(r,R,p,d)$ denote the lowest eigenvalue of the $p$--Laplacian on 
$A$, with Neumann
boundary conditions imposed where $\left|{{\bf x}}\right| = R$ and
Dirichlet
boundary conditions imposed where $\left|{{\bf x}}\right| = r$.  

Observe that 
$\alpha(r,R,p,d)$ and $\beta(r,R,p,d)$ can be written in terms of explicit Bessel functions, and 
that for any $p > 1$ $\alpha$ and $\beta$ are determined 
by ordinary differential equations, as the fundamental 
eigenfunctions are radial.  The eigenvalues $\alpha(r,R,2)$ and $\beta(r,R,2)$  and their
associated eigenfunctions are thus numerically accessible.  The Restricted Monotonicity Principle II.2 
immediately implies:
\medskip
\par\noindent
{\bf Corollary II.9 -- Annulus Bound}.
\item{a)}
Let $p>1$ and
suppose that the origin is exterior to $\Omega$, that 
$\partial \Omega_N$ is starlike with respect to the origin, and that the interior of $\Omega$ adjacent to 
$\partial \Omega_N$ is away from the origin.  Suppose further that 
$\Omega \subset \left\{{\bf x}: r < |{\bf x}| < R\right\}$.  Then 
$\lambda_1\left({\Omega, \partial\Omega,p}\right) \geq \alpha(r,R,p,d)$.

\item{b)}
Let $p>1$ and
suppose that the origin is exterior to $\Omega$, that 
$\partial \Omega_N$ is starlike with respect to the origin, and that the exterior of $\Omega$ adjacent to 
$\partial \Omega_N$ is away from the origin.  Suppose further that 
$\Omega \subset \left\{{\bf x}: r < |{\bf x}| < R\right\}$.  Then 
$\lambda_1\left({\Omega, \partial\Omega,p}\right) \geq \beta(r,R,p,d)$.

\medskip
\par\noindent
{\bf Remark}.  The Annulus Bound, as well as the Box Bound II.3, 
can be combined with Lemma II.5
to produce lower bounds analogous to (2.2) and (2.4), taking the Dirichlet boundary into 
account.

Finally, let us briefly consider the effect of changing a portion of the
boundary from $\subset \partial\Omega_N$ to $\subset \partial\Omega_D$ or {\it vice versa}.
To focus the analysis, it will be supposed that $\Omega$ is convex, and that the
unperturbed problem has boundary of one type, either
$\partial\Omega = \partial\Omega_N$ or
$\partial\Omega = \partial\Omega_D$.

The former case is easier to analyze, and indeed the 
perturbation theory when additional Dirichlet conditions 
on small subsets of positive capacity has been 
studied systematically,
for example in [Oza81], [Oza82], [Cou95], [McG96], [McG98].
A rough non--asymptotic estimate follows from Lemma II.5:
\medskip
\par\noindent
{\bf
Corollary II.10}.
{\it
Let $d \ge 2$, $1 < p < \infty$, and
suppose that $\Omega$ is strictly convex.  
Define 
$\epsilon(\Omega,p)$
as the minimum for 
${\bf x}, {\bf y} \in \partial\Omega$ and $\nu_{\bf y}$ a unit
normal to a support plane at ${\bf y}$ of ${{{\left({{\bf y} -  {\bf x}}\right) \cdot \nu_{\bf y}} \over {\left|{{\bf y} -  {\bf x}}\right|^{p+d}}}}$.
{\rm (This is a measure of the eccentricity of $\Omega$ in relation to its diameter.)}
Then}
$$\lambda_1\left({\Omega,\partial\Omega_D,p,d}\right) \ge
p^{-p}(p-1)^{p-2} {{\epsilon(\Omega,p) \left|{\partial\Omega_D}\right|}\over{\omega_d}}.$$
\medskip
\par\noindent
{\bf Proof}.
Strict convexity implies that a change of variable can be made in
the integral in 
Definition II.4 of ${1 \over {\left({m\left({\bf x}\right)}\right)^p}}$
passing from 
${\bf u} \in K \subset S^{d-1}$ to 
${\bf y} \in \partial\Omega_D$, as $K$ is the projection to the unit sphere
of $\partial\Omega_D$.  This requires a Jacobian factor 
${{{\left({{\bf y} -  {\bf x}}\right) \cdot \nu_{\bf y}} \over {\left|{{\bf y} -  {\bf x}}\right|^{d}}}},$ 
$\nu_{\bf y}$ being the outward normal where defined.  (By assumption 
in this article the outward normal is defined
a.e., but it is convenient to let $\nu_{\bf y}$ range over the normals to support planes, 
for definiteness at all boundary points and so
the formula will remain valid for some irregular domains.)
The infimum of the transformed integrand is decreased by allowing 
${\bf x}$ and ${\bf y}$ to vary independently over the closure of
$\Omega$.  In that circumstance, it is easy to see that the minimum is 
attained for 
${\bf x}, {\bf y} \in \Omega$,
yielding the result.
\hfill\qed

The problem of changing a portion of a Dirichlet boundary to Neumann
is more delicate than the converse.  
If a portion of the boundary
is changed from
having Dirichlet conditions to Neumann, then a
positive lower
bound to $\lambda_1$ can be exhibited under 
starlikeness conditions of Corollaries II.6, II.7, or II.9.  
Such bounds will, however, fail to approach 
the unperturbed eigenvalue as $\left|\partial\Omega_N\right| \rightarrow 0$.
The perturbation theory for the
introduction of Neumann 
boundaries has apparently been 
studied only under very special conditions to date
(e.g., [Oza85]).
A full theory will probably require a careful consideration of the asymptotic
behavior of eigenfunctions near the perturbation.

 \smallskip

\medskip
\noindent

\medskip 
\noindent{\bf {References}} 
\medskip

\parindent=40pt

\item{[Arr95]}
J.M. Arrieta,
 Neumann eigenvalue problems on exterior perturbations of the  domain.  
{\it J. Differential Equations} {\bf 118}  (1995)54--103.
 
 \item{[Bog07]} T. Boggio,  Sull'equazione del moto vibratorio delle  membrane elastiche, {\it Accad. Lincei, sci. fis.}, ser. 5a{\bf 16}(1907)386--393. 
 
 \item{[BuDa02]}
V.I. Burenkov and E. B. Davies,
Spectral stability of the Neumann Laplacian. (English. English summary) 
{\it J. Differential Equations} {\bf 186} (2002), no. 2, 485--508. 

\item{[CoHi37]}  
R. Courant and D. Hilbert, 
Methods of mathematical physics, vol. II, Interscience Publ., Wiley, New York, 1962.  
(Original publication 1937.)

\item{[Cou95]}
G. Courtois,
Spectrum of manifolds with holes,
{\it J. Funct. Analysis} {\bf 134}(1995)194--221.

\item{[Dav95]} E.B. Davies, {\it Spectral theory and differential operators}, Cambridge Studies
in Advanced Mathematics {\bf 42}.  Cambridge:  At the University Press, 1995.

\item{[Dav89]} E.B. Davies, {\it Heat kernels and spectral theory}, Cambridge Tracts in Mathematics
{\bf 92}.  Cambridge:  At the University Press, 1989.

\item{[DKN96]} P. Dr\'abek, A. Kufner , and F. Nicolosi,  {\it Nonlinear differential equations, Singular and degenerate case}, Pilsen, Czech Republic:  University of West Bohemia, 1996.

\item{[DKT99]} P. Dr\'abek, P. Krej\v{c}\'\i , and P. Tak\'a\v{c},  {\it Nonlinear differential equations}, Boca Raton, Florida:  CRC Press, 1999.

\item{[EdEv87]} D.E. Edmunds and W.D. Evans, {\it Spectral theory and differential operators}, Oxford:  Clarendon Press, 1987.

\item{[Ev98]}
L.C. Evans, Partial differential equations, Graduate Studies in Mathematics {\bf 19}.  Providence:  American Mathematical Society, 1998.

\item{[EvHa89]}
W.D. Evans, and D.J and Harris, On the approximation numbers of Sobolev embeddings for irregular domains, {\it Quarterly J Math Oxford} (2) {\bf 40}(1989)13--42.

\item{[FHT99]}  F. de Th\'elin, J. Fleckinger, and E.M. Harrell, II,
Boundary behavior and $L^q$ estimates for solutions of equations containing the 
$p$--Laplacian, with J. Fleckinger and F. de Th\'elin, 
{\it Electronic J. Diff. Eqns.} {\bf 1999}(1999)1--19.

\item{[Fra79]}  L.E. Fraenkel, On regularity of the boundary in the theory of Sobolev 
spaces.  {\it Proc. Lond. Math. Soc.} (3) {\bf 39}(1979)385--427.

\item{[Har20]}  G.H. Hardy, Note on a theorem of Hilbert.  {\it Math. Z.} {\bf 6}(1920)314--317.

\item{[HLP34]} G.H. Hardy, J.E. Littlewood, and G. P\'olya, {\it Inequalities}.  Cambridge:  At the University Press, 1959.  Original publication 1934.

\item{[Ha93]}
E.M. Harrell II,
Lecture at the Workshop on Partial Differential
Equations and Fractals, Toulouse, France, May, 1993, unpublished.

\item{[HSS91]}
R. Hempel, L.A. Seco, and B. Simon,
The essential spectrum of Neumann Laplacians on some bounded singular domains, 
{\it J. Funct. Analysis} {\bf 102}(1991)448--483.

\item{[Ken94]}
C. Kenig, 
{\it Harmonic analysis techniques for second order elliptic boundary value problems}, 
Conference Board of the Mathematical Sciences, Regional Conference Series in Mathematics 83.  Providence:  American Mathematical Society, 1994.

\item{[McG96]}
I. McGillivray, Capacitary estimates for Dirichlet eigenvalues, 
{\it J. Funct. Analysis} {\bf 139}(1996)244--259.

\item{[McG98]}
I. McGillivray, Capacitary asymptotic expansion of the groundstate to second order, Commun. P.D.E. {\bf 23}
(1998)2219--2252

\item{[Maz81]}  
V.G. Maz'ja, 
Sobolev spaces.  New York:  Springer, 1985.  First published 1981.

\item{[Maz73]}  
V.G. Maz'ja, 
On $(p,l)$--capacity, embedding theorems, and the spectrum of a self-adjoint elliptic operator, 
{\it Izv. Akad. Nauk SSSR ser. Mat.} {\bf 37}(1973)356-385:  {\it Transl. Math. USSR IZV} 7(1973)357-387.

\item{[OpKu90]} 
B. Opic and A. Kufner, 
{\it Hardy--type inequalities,} Pitman Research Notes in Math. {\bf 219} Boston: Longman, 1990.

\item{[Oza81]}
S. Ozawa, 
Singular variation of domains and eigenvalues of the Laplacian.  
{\it Duke Math. J.} {\bf 48}(1981)767--778.

\item{[Oza82]}
S. Ozawa, An Asymptotic Formula for the Eigenvalues of the Laplacian in a Domain with a Small Hole, {\it Proc. Japan Acad.} {\bf 58}, ser. A (1982)5--8.

\item{[Oza85]}
S. Ozawa
Shin
Asymptotic property of an eigenfunction of the Laplacian under singular  variation of domains--the Neumann condition.  
{\it Osaka J. Math.} {\bf 22}  (1985)639--655.
%
%
%
%
%
%
\end